\documentclass [a4paper,twoside,12pt]{article}

\usepackage[utf8]{inputenc}
\usepackage{authblk}

\usepackage{amsmath,amsfonts,amssymb}
\usepackage{vmargin,graphicx,theorem}

\usepackage[english]{babel}
\usepackage{enumerate}
\usepackage{color}
\usepackage{pst-fill,pst-grad,pst-plot,pst-eucl,pstricks-add,pst-node}
\RequirePackage[colorlinks,linkcolor=blue,citecolor=blue,urlcolor=blue]{hyperref}

\setpapersize[portrait]{A4}
\setmarginsrb{2cm}{0.7cm}{1.5cm}{2cm}{1.5cm}{0cm}{0.5cm}{2cm}

\newcommand{\dR}{\ensuremath{\mathbf{R}}}

\def\Hess{\mathop{\rm Hess}\nolimits}

\newtheorem{ethm}{Theorem}[section]

\newtheorem{erem}[ethm]{Remark}

\newcommand{\proofend}{~$\rhd$}
\newcommand{\proofbegin}{~$\lhd$}

\newtheorem{thmx}{Theorem}

\newcommand{\PAR}[1]{\ensuremath{{\left(#1\right)}}} 
\newcommand{\SBRA}[1]{\ensuremath{{\left[#1\right]}}} 

\renewcommand{\phi}{\varphi}

\renewcommand{\geq}{\geqslant}

\newcommand{\N}{\ensuremath{\mathbf{N}}}
\newcommand{\R}{\dR}

\newcommand{\grad}{{\rm grad}}

\renewcommand{\S}{\ensuremath{\mathbb{S}}}

\newcommand{\beq}{\begin{equation}}\newcommand{\eeq}{\end{equation}}
\begin{document}

	\title{Sobolev's inequality under a curvature-dimension condition}
	
	\author{Louis Dupaigne, Ivan Gentil, Simon Zugmeyer}
	
	\date{\today}
	\maketitle

	\begin{center}
		Résumé
	\end{center}
	
	Dans cette note, nous proposons une nouvelle preuve de l'inégalité de Sobolev sur les variétés à courbure de Ricci minorée par une constante positive. Le résultat avait été obtenu en 1983 par Ilias. Nous présentons une preuve très courte de ce théorème, dressons l'état de l'art pour cette fameuse inégalité et expliquons en quoi notre méthode, qui repose sur un flot de gradient, est simple et robuste. En particulier, nous élucidons les calculs utilisés dans des travaux précédents, à commencer par un célèbre article de Bidaut-Véron et Véron publié en 1991.

	\begin{center}
		Abstract
	\end{center}
	In this note we present a new proof of Sobolev's inequality under a uniform lower bound of the Ricci curvature. This result was initially obtained in 1983 by Ilias. Our goal is to present a very short proof, to give a review of the famous inequality and to explain how our method, relying on a gradient-flow interpretation, is simple and robust.  In particular, we elucidate computations used in numerous previous works, starting with Bidaut-V\'eron and V\'eron's 1991 classical work.

	\setmarginsrb{2cm}{0.7cm}{1.5cm}{2cm}{1.5cm}{0cm}{0.5cm}{2cm}
	
	\section{Introduction}
	
	Given $d\in\N$, $d\ge2$, and $p\in [1,d)$, let $p^*\in[1,+\infty)$ denote Sobolev's exponent, that is
	$$
	\frac1{p^*} = \frac1p -\frac 1d.
	$$
	According to Sobolev's inequality, there exists a constant $A>0$ such that for every $\varphi\in C^\infty_c(\R^d)$, 
	$$
	\Vert \varphi\Vert_{L^{p^*}(\R^d)}\le A \Vert \nabla \varphi\Vert_{L^{p}(\R^d)},
	$$
	see~\cite{Sob38}, as well as~\cite{Gag58, Nir59} for the case $p=1$,~\cite{Rod66, Aub76, Tal76} for the value of the sharp constant $A$ and the expression of the extremals,~\cite{Lie83} for a more direct proof using rearrangements and~\cite{CafGidSpr89} for the classification of all positive solutions to the associated Euler-Lagrange equation. 
	In the special case $p=2$, using the stereographic projection (see e.g.~\cite{LeePar87}), the sharp Sobolev inequality in $\R^d$ is equivalent to 
	\begin{equation*}
	\frac1{q-2}\left( ||v||_{L^{q}(\S^d)}^2 - ||v||_{L^2(\S^d)}^2\right)\le\frac{1}{d} ||\nabla v||_{L^2(\S^d)}^2,
	\end{equation*}
	where $q=2^*$, $v\in C^\infty(\S^d)$, $\S^d$ is the standard sphere equipped with its normalized\footnote{in other words, the normalized measure $\nu$ is proportional to the Riemannian volume and $\nu(\S^d)=1$} measure. The inequality is again sharp and the extremals are known, see~\cite{Aub76b}, as well as Theorem 5.1 p.~121 in~\cite{Heb99}. In fact, the inequality is true for every $q\in[1,2^*]$, $q\neq2$, see~\cite{BidVer91, Bec93}, as well as~\cite{Dem05}, Section~3.11 for the case $q\in[1,2)$. Also note that letting $q\to2$, one recovers the sharp log-Sobolev inequality. In~\cite{Ili83}, Sobolev's inequality has been generalized as follows to any compact Riemannian manifold $(M,\mathfrak g)$ with positive Ricci curvature. 
	
	\begin{thmx}[\cite{Ili83}]
		\label{th:sobolev}
		Let $(M,\mathfrak g)$ be a smooth connected, compact, $d$-dimensional Riemannian manifold, $d\ge3$. Assume that the Ricci curvature of $M$ is uniformly bounded from below by a constant $\rho>0$. Let $q=2^*=\frac{2d}{d-2}$. Then, for all $v\in C^\infty(M)$,
		\begin{equation}\label{sobolev:M}
		\frac1{q-2}\left( ||v||_{L^{q}(M)}^2 - ||v||_{L^2(M)}^2\right)\le\frac{1}{d}\frac{d-1}{\rho} ||\nabla v||_{L^2(M)}^2,
		\end{equation}
		where $M$ is equipped with its normalized measure.
	\end{thmx}
	
	\begin{erem}
		It is not necessary to assume that $M$ is compact, as follows from Myer's theorem (see e.g.~\cite{Heb99} p.~100 for a geometric proof, or combine Theorems 3.2.7, 6.6.1 and 6.8.1 in~\cite{bgl-book} for an analytic proof).
	\end{erem}
	
	Many proofs of Theorem~\ref{th:sobolev} are available. The approach in~\cite{Ili83} relies on symmetrization arguments and the Lévy-Gromov isoperimetric inequality~\cite{Gro07}, the rigorous proof of which seems involved, see e.g.~\cite{Vil19}. 
	The proof of~\cite{BidVer91} clarifies computations of~\cite{GidSpr81}, but does not elucidate them. The latter paper presumably took inspiration from Obata's work~\cite{Oba62} (also described in~\cite{BerGauMaz71}, pp.~179--185). 
	In~\cite{BakLed96} (see Theorem 6.10 p.~107 in~\cite{bakry94} for the actual proof, as well as Chapter~6 in~\cite{bgl-book} for a more recent and thorough account), the inequality is generalized to any Markov generator satisfying the curvature-dimension condition $CD(\rho,n)$, $\rho>0$, $n>2$. Among other tools, their proof makes use of the Bakry-\'Emery method (or $\Gamma$-calculus) and a rather unintuitive change of unknown which was already present in the aforementioned litterature. The proof of Fontenas~\cite{Fon97} provides a sharper version of the inequality in terms of the generator's best Poincar\'e constant in the case $q\in[2,2^*)$. His computations, inspired by~\cite{Rot86}, use again the $\Gamma$-formalism and recast the proof in a yet simpler form, but still fail short of making it transparent. 
	In~\cite{DelDol02}, Sobolev's inequality in $\R^d$ appears as a limiting case of a family of optimal Gagliardo-Nirenberg inequalities. This paper puts forward two important tools for our purposes: the classification of solutions to the associated Euler-Lagrange, based here on the symmetry result of~\cite{GidNiNir81} and, more importantly, the connection between Sobolev's inequality and the convergence to equilibrium of solutions to the fast-diffusion equation, or rather to a Fokker-Planck-type equation obtained by rescaling. 
	The fast-diffusion and porous medium equations had just been reformulated in~\cite{otto2001} as a gradient flow in Wasserstein space, leading the way to the reinterpretation of Sobolev's inequality (and more generally the Gagliardo-Nirenberg inequalities studied by del Pino and Dolbeault) as a simple convexity inequality along a flow, in other words as an entropy-entropy production inequality. This latter point of view was taken in~\cite{carTos00},~\cite{cJM+01} and~\cite{carVaz03} to establish  Sobolev-type inequalities in $\R^d$ and more recently simplified and generalized to convex euclidean domains  in~\cite{Zug19}. 
	Soon after,~\cite{CorNazVil04} gave a short proof using optimal transport, but valid in the euclidean setting only. 
	The extension of the Bakry-\'Emery method to nonlinear flows was further cleverly extended in the Riemannian setting in~\cite{demange2008}, although without Otto's geometric insight, but with a twist: the use of two distinct entropy functionals, the evolutions of which can be related through a simple differential inequality. Other recent generalizations include the cases of $\text{\rm RCD}^*(\rho,n)$-spaces~\cite{Pro15},  $\text{\rm CD}^*(\rho,n)$-spaces~\cite{cavmon} and Riemannian manifolds with boundary~\cite{IliSho18}.
	Going back to the euclidean setting, but with weights,~\cite{DolEstLos16} extended the method to prove the sharp Caffarelli-Kohn-Nirenberg inequalities and the associated Liouville-type results.
	
	We probably forgot to cite important contributions and judging by the extent of the bibliography, one may wonder why we intend to give here yet another proof of Sobolev's inequality. From our point of view, the proof presented below, inspired by~\cite{DolEstLos16}, has the advantage of being short, transparent and hopefully robust. In particular, with no extra work, our proof yields the following generalization\footnote{For convenience of the reader, in Section~\ref{proofs}, we recall the definition of the $CD(\rho,n)$ condition used in Theorem~\ref{th:B}.} of Theorem~\ref{th:sobolev}, due to~\cite{BakLed96}.
	
	\begin{thmx}[\cite{BakLed96}]\label{th:B}
		Let $\alpha\in(0,1)$. Assume that $(M,\mathfrak g)$ is a $C^{2,\alpha}$, compact, connected, $d$-dimensional Riemannian manifold, $d\ge1$. Let $W\in C^{2}(M;\R)$ and $L=\Delta - \nabla W\cdot\nabla$ satisfy the $CD(\rho,n)$ condition for some $\rho>0$ and $n\in[d,+\infty)$, $n>2$. Let $q=\frac{2n}{n-2}$. Then, for all $v\in C^\infty(M)$,
		\begin{equation*}
		\frac1{q-2}\left( ||v||_{L^{q}(M)}^2 - ||v||_{L^2(M)}^2\right)\le\frac{1}{n}\frac{n-1}{\rho} ||\nabla v||_{L^2(M)}^2,
		\end{equation*}
		where $M$ is equipped with the measure $d\nu= \frac{e^{-W}}Z d\text{\rm Vol}_g$, with $Z\in \R_+^*$ chosen so that $\nu(M)=1$.
	\end{thmx}
	
	\begin{erem}
		Again, it is not necessary to assume that $M$ is compact, as follows from the generalized Myer's theorem proved in~\cite{BakLed96}. 
	\end{erem}
	As another by-product of our proof, we obtain the following rigidity result, improving previous results given in~\cite{GidSpr81, BidVer91, LicVer95, LicVer98, BakLed96, Fon97, DolEstKowLos14, DolEstLos14}, which, as stated, seems new.
	
	\begin{ethm}
		\label{th:liouville}
		Assume that $(M,\mathfrak g)$ is a $C^{2}$, compact, connected, $d$-dimensional Riemannian manifold, $d\ge1$. Let $W\in C^{2}(M;\R)$ and $L=\Delta - \nabla W\cdot\nabla$ satisfy the $CD(\rho,n)$ condition for some $\rho>0$ and $n\in[d,+\infty)$, $n>2$. Let $q=\frac{2n}{n-2}$. Assume that $v\in C^2(M)$, $v>0$, is a nonconstant solution to 
		\begin{equation}\label{ele}
		-A\, L v + v = v^{q-1}f(v)\quad\text{in $M$},
		\end{equation}
		where  $A>0$ and $f\in C^{1,\alpha}(\R_+^*;\R_+^*)$, $\alpha\in(0,1)$ and $f$ is nonincreasing. Let $A^*=\frac{4(n-1)}{n(n-2)\rho}$. Then, $A\le A^*$. In addition, if $A=A^*$, then $f$ is constant on $[0,\Vert v\Vert_\infty]$.
	\end{ethm}
	
	\begin{erem} If equality holds in Sobolev's inequality \eqref{sobolev:M} for some nonconstant function $v$, then $v$ solves the associated Euler-Lagrange (equation \eqref{ele} with $n=d$, $A=A^*$, $L=\Delta$ and $f$ constant). As follows from the proof of Theorem \ref{th:liouville}, the function $\Phi=v^{-\frac4{d-2}}$ solves the equation $\nabla^2\Phi = \frac{\Delta\Phi}{d}\mathfrak g$ in $M$. This in turn implies that $(M, \mathfrak g)$ is conformally diffeomorphic to the round sphere, see e.g. Lemme 6.4.3 in \cite{Heb97}. If we assume in addition that $(M, \mathfrak g)$ is Einstein, letting $d_{\mathfrak g}$ denote its Riemannian distance, we have in fact that $(M, \mathfrak g)$ is isometric to the round sphere and that $v(x)=(\beta-\cos(d_{\mathfrak g}(x_0,x))^{-\frac{d-2}2}$ for some $\beta>1$ and $x_0\in M$, see e.g. Theorem 5.1 in \cite{Heb99} and its proof.
	\end{erem}

	
	
	\section{Proofs of Theorem~\ref{th:sobolev}, Theorem~\ref{th:B} and Theorem~\ref{th:liouville}}\label{proofs}
	\subsection{Proof of Theorem~\ref{th:sobolev}}
	\label{sec-2.1}
	Fix $q\in[1,2^*)$. By the (non-sharp but tight) Sobolev inequality, there holds
	\begin{equation}\label{eq:sobno}
	\Vert v\Vert_q^2 \le A\Vert\nabla v\Vert_2^2 + \Vert v\Vert_2^2,
	\end{equation}
	for some $A\in\R_+^*$ and every $v\in H^1(M)$, apply e.g.~\cite{Heb99}, Corollary 2.1 and \cite{bgl-book}, Proposition~6.2.2. Given $A\in\R_+^*$, consider the minimization problem
	$$
	I(A)=\inf \left\{A ||\nabla v||_2^2+||v||_2^2\;:\; v\in H^1(M)\;,\; ||v||_{q}=1\right\}.
	$$
	Then,~\eqref{eq:sobno} holds if  $I(A)=1$. Thanks to the Banach-Alaoglu-Bourbaki and Rellich-Kondrakov compactness theorems (see e.g.~\cite{Bre11} Theorem 3.16 and~\cite{Heb99} Theorem 2.9), there exists a minimizer $v\in H^1(M)$ s.t. $||v||_{q}=1$. By Stampacchia's theorem~\cite{Sta66}, $\vert v\vert$ is also a minimizer, so we may assume that $v\ge 0$ a.e. in $M$. In addition, a  constant multiple of $v$ (abusively denoted the same below) is a weak solution to
	\begin{equation}
	\label{10}
	-A\Delta v + v = v^{q-1}\quad\text{in $M$}.
	\end{equation}
	By standard elliptic regularity (see e.g.~\cite{Heb97}, proof of Theorem 6.2.1, p.~248) $v\in C^3(M)$ and by the strong maximum principle (see e.g.~\cite{Heb97}, Theorem 5.7.2), $v>0$ in $M$.
	
	Define the pressure function $\Phi=v^{-\frac{q-2}2}$. Then, $\Phi$ solves
	\begin{equation}\label{eq:p}
	\Phi\Delta \Phi -\frac {d'}2\vert\nabla \Phi\vert^2 = -\lambda (\Phi^2-1)\quad\text{in $M$,}
	\end{equation}
	where $d'=\frac{2q}{q-2}$ and $\lambda=\frac{q-2}{2A}=\frac2{(d'-2)A}$.
	Multiply  equation~\eqref{eq:p} by \(\Delta \Phi^{1-{d'}}\) and integrate. For the right-hand-side we find,
	\begin{align*}
	\int {\lambda(\Phi^2-1)}\Delta \Phi^{1-d'}
	&=  \lambda\int \Phi^2\Delta \Phi^{1-d'}=-\lambda\int \nabla \Phi^2\cdot\nabla \Phi^{1-d'} \\
	&= 2\lambda(d'-1)\int \vert\nabla \Phi\vert^2 \Phi^{1-{d'}}=c\int \Gamma(\Phi)\Phi^{1-d'}
	\end{align*}
	where we expressed the carr\'e du champ operator
	$
	\Gamma(\Phi) = 
	\vert\nabla \Phi\vert^2
	$ 
	and where $c=2\lambda(d'-1)=4\frac{d'-1}{(d'-2)A}$.
	For the left-hand side, we obtain
	\begin{align*}
	\int \PAR{\Phi\Delta \Phi - \frac {d'}2\vert{\nabla \Phi}\vert^2}\Delta \Phi^{1-d'}
	&=  \int \Delta\PAR{\Phi\Delta \Phi - \frac {d'}2\vert{\nabla \Phi}\vert^2}\Phi^{1-d'} \\
	&= \int \left[{(\Delta \Phi)^2 + \Phi\Delta^2\Phi + 2\nabla \Phi\cdot\nabla\Delta \Phi}- \frac {d'}2{\Delta\vert{\nabla \Phi}\vert^2}\right]\Phi^{1-{d'}}\\
	&= - \int \PAR{{d'}\Gamma_2(\Phi)-(\Delta \Phi)^2}\Phi^{1-{d'}},
	\end{align*}
	where we expressed the iterated carr\'e du champ 
	$
	\Gamma_2(\Phi) =\frac12\Delta \vert\nabla \Phi\vert^2 - \nabla \Phi\cdot\nabla\Delta \Phi
	$
	and used the fact that
	\[
	\int \Phi^{2-{d'}}\Delta^2 \Phi = ({d'}-2) \int (\nabla \Phi\cdot\nabla\Delta \Phi) \Phi^{1-{d'}}.
	\]
	Collecting the left and right-hand sides and dividing by ${d'}$, we find
	\begin{equation}
	\label{eq:gamma2}
	\int\PAR{\Gamma_2(\Phi)-\frac 1{d'}(\Delta \Phi)^2-\frac c{d'}\Gamma(\Phi)}\Phi^{1-{d'}} =0.
	\end{equation}
	The celebrated Bochner-Lichnerowicz formula 
	states\footnote{and motivates the definition of $\Gamma_2$} that
	$$
	\Gamma_2(\Phi) = \Vert  \nabla^2\, \Phi\Vert_{H.S}^2+\text{\rm Ric}_{\mathfrak g}(\nabla \Phi, \nabla \Phi),
	$$
	where  $\nabla^2 \Phi$ denotes the Hessian of $\Phi$,  $\Vert  \nabla^2\, \Phi\Vert_{H.S}^2$ the square of its Hilbert-Schmidt norm (the sum of the squares of its components) and $\text{\rm Ric}_{\mathfrak g}$ the Ricci tensor of the Riemannian manifold $(M,\mathfrak g)$.
	Using the Cauchy-Schwarz inequality on the one hand and the assumption $\text{\rm Ric}\ge\rho {\mathfrak g}$ on the other hand, we find
	$$
	\Gamma_2(\Phi)\ge \frac1d(\Delta \Phi)^2 + \rho\Gamma(\Phi) 
	$$
	and so
	$$
	\left(\frac 1d -\frac 1{d'}\right)  \int(\Delta \Phi)^2\Phi^{1-{d'}} + \left(\rho-\frac c{d'}\right)\int \Gamma(\Phi)\Phi^{1-{d'}}\le0.
	$$
	Since $q< 2^*$, we have $d< {d'}$ and so, if  $\rho\ge\frac c{d'}$ i.e.
	$$
	A\ge \frac{4({d'}-1)}{{d'}({d'}-2)\rho},
	$$
	we deduce that $\Delta \Phi =0$ in $M$. Integrating against $\Phi$, $\Phi$ is constant. Hence $v=1$, $I(A)=1$, and~\eqref{eq:sobno} holds for  $A=\frac{4({d'}-1)}{{d'}({d'}-2)\rho}$.
	Let $q\nearrow 2^*$. Then ${d'}\searrow d$ and~\eqref{sobolev:M} follows.
	
	\subsection{Proof of Theorem~\ref{th:B}}
	\subsubsection{The $CD(\rho,n)$ condition.}
	Let us quickly explain the definitions and notations used in the theorem. Clearly, a second order differential operator of the form\footnote{Here $\Delta$ is the  Laplace-Beltrami operator on $(M,\mathfrak g)$, the dot product designates the Riemannian metric $\mathfrak g$ and $\vert\cdot\vert$ the associated norm.}
	$L=\Delta - \nabla W\cdot\nabla$ fails to satisfy the chain rule: if $\Phi\in C^2(M)$ is not constant, $L(\Phi^2)\neq 2\Phi L\Phi$. The defect is measured by the carr\'e du champ operator defined for $\Phi\in C^2(M)$ by
	$$
	\Gamma(\Phi) = \frac12 L(\Phi^2) -\Phi L\Phi.
	$$
	By a simple and direct computation, $\Gamma(\Phi)=\vert\nabla \Phi\vert^2$. Abusing notation slightly, we let $\Gamma(\Phi,\Psi)=\nabla \Phi\cdot\nabla \Psi$ denote the polar form of $\Gamma$. Now, repeat the above consideration by replacing the product of real numbers, seen as a bilinear form, by the carré du champ operator $\Gamma$: again $L$ fails to satisfy the chain rule and the defect is measured by the iterated carré du champ operator, defined for $\Phi\in C^3(M)$ by
	\begin{equation}
	\label{eq:gamma22}
	\Gamma_2(\Phi) = \frac12  L(\Gamma(\Phi))-\Gamma(\Phi,L\Phi).
	\end{equation}
	Thanks to the  Bochner-Lichnerowicz formula, the $\Gamma_2$ operator can be computed as follows:
	$$
	\Gamma_2(\Phi) =  \Vert\nabla^2 \Phi\Vert_{H.S.}^2 + (\text{\rm Ric}_{\mathfrak g}+ \nabla^2 W)(\nabla \Phi, \nabla \Phi) .
	$$
	Given, $\rho\in\R$ and $n\in[d,+\infty]$, the operator $L$ is then said to satisfy the $CD(\rho,n)$ condition if for every $\Phi\in C^3(M)$,
	\begin{equation}
	\label{100}
	\Gamma_2(\Phi)\ge \rho\Gamma(\Phi)+\frac1n(L\Phi)^2.
	\end{equation}
	Note that when $W=0$, $L\Phi=\Delta \Phi$. By the Cauchy-Schwarz inequality\footnote{with equality if and only if $\nabla^2 \Phi = \frac{\Delta \Phi}{d}\mathfrak g$.}, $\Vert\nabla^2 \Phi\Vert_{H.S.}^2 \ge \frac1d(\Delta \Phi)^2$ so that, in this case, the $CD(\rho,d)$ condition\footnote{We recall that $d$ is the dimension of $M$} is equivalent to the lower bound $\text{\rm Ric}_{\mathfrak g}\ge \rho {\mathfrak g}$. 
	\subsubsection{Proof of Theorem~\ref{th:B}}
	Let us review the proof of Theorem~\ref{th:sobolev}. We start similarly with the tight but non-sharp Sobolev's inequality~\eqref{eq:sobno}, the proof of which remains unchanged (e.g. adapt~\cite{Heb99} Theorem 4.1). Since $M$ is compact and $W$ continuous, ${e^{-W}}$ is bounded above and below by positive constants. So, the Riemannian volume and the measure $d\nu= \frac{e^{-W}}Z d\text{\rm Vol}_g$ yield the same Sobolev space $H^1(M,d\nu)=H^1(M,d\text{\rm Vol}_g)$. In particular, by the same proof, the quantity $I(A)$ has a nonnegative minimizer $u$, which this time solves
	$$
	-A\, L v + v = v^{q-1}\quad\text{in $M$},
	$$
	leading to
	\begin{equation*}
	\Phi L \Phi -\frac{n'}2\vert\nabla \Phi\vert^2 = -\lambda (\Phi^2-1)\quad\text{in $M$},
	\end{equation*}
	where the definition of $\Phi$ is unchanged, $n'=\frac{2q}{q-2}$ and $\lambda=\frac{q-2}{2A}=\frac2{(n'-2)A}$. Multiply by $L(\Phi^{1-n'})$ and integrate. Using the formulas $\int_M (Lu) v\;d\nu=\int_M uLv\;d\nu=-\int_M \Gamma(u,v)\;d\nu$, the exact same computations lead to  
	\begin{equation*}
	\label{eq:gamma}
	\int\PAR{\Gamma_2(\Phi)-\frac 1{n'}(L \Phi)^2-\frac c{n'}\Gamma(\Phi)}\Phi^{1-n'}d\nu =0,
	\end{equation*}
	where $c=2\lambda(n'-1)=4\frac{n'-1}{(n'-2)A}$. Now apply the $CD(\rho,n)$ condition to deduce that~\eqref{eq:sobno} holds for  $A=\frac{4(n'-1)}{n'(n'-2)\rho}$.
	Let $q\nearrow \frac{2n}{n-2}$. Then, $n'\searrow n$ and the theorem follows.

	\subsection{Proof of Theorem~\ref{th:liouville}}
	Repeating once again the above computation we arrive at
	$$
	\int\left(\Gamma_2(\Phi)-\rho\Gamma(\Phi)-\frac1n(L\Phi)^2\right)\;d\nu +
	\left(\rho-\frac cn\right)\int \Gamma(\Phi)\Phi^{1-n}d\nu+\lambda\int f'(v)\Phi^2\nabla v\cdot\nabla{\Phi^{1-n}}d\nu=0,
	$$	 
	where $c=2\lambda(n-1)=4\frac{n-1}{(n-2)A}$ and $\lambda=\frac{q-2}{2A}=\frac2{(n-2)A}$.
	By the $CD(\rho,n)$ condition, the first integral is nonnegative. Since $f$ is nonincreasing, the last integral is also nonnegative. Finally, the coefficient in front of the second integral is strictly positive if $A>A^*$, so that $v$ must be constant in that case. If $A=A^*$, then all the first and third integrals vanish. In particular, $f$ is constant on $[0,\Vert v\Vert_\infty]$.
	
	%
	%
	
	\section{Sobolev's inequality is a convexity inequality for Renyi entropies in Wasserstein space}
	\label{sec-3}
	
	In this section,  we explain the genesis of our short proof of Theorems~\ref{th:sobolev} and~\ref{th:liouville}. Our strategy consists in using a gradient flow defined on the set of probability measures over $M$, equipped with the Wasserstein distance. If one uses the appropriate functionals, the proof is rather simple. In the next paragraph, we explain first how a gradient flow in the usual Euclidean space $\R^m$ can be used to derive sharp convexity inequalities. The extension of the method to the Wasserstein space is next presented in Section~\ref{sec-3.2}. The computations are not new, but this presentation and this point of view seem to be new and useful. 
	
	Some of our considerations will be formal: although this can be done, we do not try to make all arguments rigorous, but we provide references to do so.  Instead, we ask the reader to keep in mind that we only want to give a guideline to the rigorous proofs presented previously.

	
	\subsection{A review of gradient flows in Euclidean space}
	\label{sec-3.0}

	Let $m\ge1$ and $F:\R^m\mapsto\R$ any $C^2$ function, that we call {\it entropy} in what follows. Assume that $F$ is strictly convex and coercive i.e. 
	$
	\lim_{\vert x\vert\to+\infty} F(x) = +\infty.
	$
	Then, $F$ has unique critical point $x^*$. In addition,
	$$
	F(x^*)=\inf_{x\in\R^m}F(x).
	$$
	In order to locate the point of minimum $x^*$, one can start from an arbitrary point $x\in\R^m$ and follow the gradient flow associated to $F$. More precisely, let $t\mapsto S_t(x)$ denote the solution of the ODE
	\begin{equation}
	\label{213}
	\left\{
	\begin{array}{l}
	\displaystyle\frac{d}{dt}S_t(x)=-\nabla F(S_t(x))\\
	\displaystyle S_0(x)=x.
	\end{array}
	\right.
	\end{equation}
	Thanks to the Cauchy-Lipschitz theorem, $t\mapsto S_t(x)$ is well-defined on a maximal interval $I$ containing $t=0$. In fact, the solution is bounded, hence global, since $F$ is coercive and nondecreasing along the flow:
	\begin{equation}\label{lyap}
	\frac d{dt}F(S_t(x)) = -\vert \nabla F(S_t(x))\vert^2\le 0.
	\end{equation}
	In addition, given any $x\in\R^m$, 
	\begin{equation}
	\label{214}
	\lim_{t\rightarrow\infty}S_t(x)=x^*.
	\end{equation}
	Indeed, since $F$ is bounded below and \eqref{lyap} holds, there exists a sequence $t_n\to+\infty$ such that $\vert \nabla F(S_{t_n}(x))\vert\to0$. Since $(S_{t}(x))$ is bounded, up to extraction, $(S_{t_n}(x))$ also converges and by continuity of $\vert\nabla F\vert$, its limit must be $x^*$. Using \eqref{lyap} once more, $F(S_t(x))\le F(S_{t_n}(x))$ for $t\ge t_n$ and so $F(S_t(x))$ decreases to $F(x^*)$. \eqref{214} follows.
	
	If we further assume that $F$ is strongly convex, i.e. $\nabla^2 F\ge \rho\,\text{\rm Id}$ for some $\rho>0$, then the rate of convergence of the {\it entropy along its gradient flow} can be quantified (as we shall prove shortly):
	$$
	F(S_t(x)) - F(x^*) \le e^{-2\rho t}\left(F(x) - F(x^*)\right).
	$$
	Note that equality holds when $t=0$ and so we can differentiate the inequality at $t=0$. This yields the following equivalent convexity inequality
	$$
	F(x) - F(x^*)\le \frac1{2\rho}\vert\nabla F(x)\vert^2.
	$$
	Note that the inequality is sharp in the sense that it is an equality for $F(x)=\rho\vert x\vert^2/2$.
	In fact, one can be a bit more general and consider the following convexity inequality
	\begin{equation}
	\label{216}
	G(x^*)\leq \frac{1}{2\rho}|\nabla F(x)|^2+G(x),
	\end{equation}
	which holds true whenever $G\in C^2(\R^m)$ and $F$ satisfy the following convex condition: there exits $\rho>0$ such that uniformly in $\R^m$, 
	\begin{equation}
	\label{215}
	\nabla F\cdot \nabla^2 F\,\nabla F\geq -\rho\nabla F\cdot \nabla G.
	\end{equation}
	We provide three proofs of this fact, ending with the most robust.

	\begin{enumerate}[\bf 1.]
		\item {\bf  A direct proof based on the gradient flow. } Differentiating ~\eqref{213} once more, gives, for any $x\in\R^m$, 
		\begin{multline*}
		\frac{d^2}{dt^2}F(S_t(x))=2\nabla F(S_t(x))\cdot\nabla^2 F(S_t(x))\nabla F(S_t(x))\geq\\
		-2\rho\nabla F(S_t(x))\cdot\nabla G(S_t(x))=2\rho\frac{d}{dt}G(S_t(x)).
		\end{multline*}
		Integrating over $[0,\infty]$ the previous inequality becomes, 
		$$
		\int_0^\infty\frac{d^2}{dt^2}F(S_t(x))dt\geq 2\rho\int_0^\infty\frac{d}{dt}G(S_t(x))dt.
		$$
		Since
		\begin{equation}
		\label{217}
		\lim_{t\rightarrow\infty} \vert\nabla F (S_t(x))\vert=0,
		\end{equation}
		we have 
		$$
		-\frac{d}{dt}F(S_t(x))\big|_{t=0}\geq2\rho(G(x^*)-G(x)).
		$$
		Since $-\frac{d}{dt}F(S_t(x))|_{t=0}=|\nabla F(x)|^2$, we proved the inequality~\eqref{216}, under the condition~\eqref{215}.
		
		As we can see, inequality~\eqref{216} is just a clever convex inequality under the convex condition~\eqref{215}. As we shall see, when generalizing this proof to an infinite-dimensional setting, we are faced with two problems: proving rigorously the existence of the gradient flow $(S_t)_{t\geq0}$ and proving the two limits~\eqref{214} and~\eqref{217}. 
		
		\item {\bf A proof based on a minimization problem and the gradient flow.}\label{3.2} To prove~\eqref{216}, we fix a constant $A>0$, compute the quantity 
		$$
		I(A):=\inf_{x\in\R^m}\SBRA{A |\nabla F(x)|^2+G(x)}
		$$
		and show that  for $A>\frac1{2\rho}$, $G(x^*)\le I(A)$. Letting $A\searrow\frac1{2\rho}$, \eqref{216} will then follow. If $G$ is coercive, which we assume in this approach, then there exits $\bar x\in\R^m$  such that 
		\begin{equation}
		\label{227}
		\inf_{x\in\R^m}\SBRA{A|\nabla F(x)|^2+G(x)}=A|\nabla F(\bar x)|^2+G(\bar x).
		\end{equation}
		We now consider $(S_t(\bar x))_{t\geq0}$, the gradient flow starting from $\bar x$.  Then, since $\bar x$ is a minimizer, we have 
		$$
		\frac{d}{dt}\SBRA{A|\nabla F(S_t(\bar x))|^2+G(S_t(\bar x))}\Big|_{t=0}\geq0. 
		$$
		In addition, 
		\begin{multline}
		\label{218}
		\frac{d}{dt}\SBRA{A |\nabla F(S_t(\bar x))|^2+G(S_t(\bar x))}\Big|_{t=0}=-2A\nabla F(\bar x)\cdot\nabla^2F(\bar x)\nabla F(\bar x)-\nabla G(\bar x)\cdot \nabla F(\bar x)=\\
		\left[-\frac1\rho \nabla F(\bar x)\cdot\nabla^2F(\bar x)\nabla F(\bar x)-\nabla G(\bar x)\cdot\nabla F(\bar x)\right]-\left(2A-\frac1{\rho}\right)\nabla F(\bar x)\cdot\nabla^2F(\bar x)\nabla F(\bar x).
		\end{multline}
		
		Since $F$ is strictly convex and \eqref{215} holds, if $A>\frac1{2\rho}$, we see that if
		$\bar x\neq x^*$, $\nabla F(x)\neq0$ and so
		$$
		\frac{d}{dt}\SBRA{A\,|\nabla F(S_t(\bar x))|^2+G(S_t(\bar x))}\Big|_{t=0}<0, 
		$$
		which is impossible since $\bar x$ is a minimizer.  Hence,  $\bar x=x^*$ and the following inequality holds, 
		$$
		G(x^*)\leq A |\nabla F(x)|^2+G(x),
		$$
		for any $x\in\R^d$ and $A>\frac1{2\rho}$. This proves the desired inequality  \eqref{216}, by letting $A\to \frac1{2\rho}$. 
		Note that in this approach, we no longer need to prove the asymptotic behavior of the gradient flow $(S_t)_{t\geq0}$  but we still need to know its existence. 
		\item {\bf A proof based  on the minimization problem only. } As in the previous proof, let $\bar x$ given by equation~\eqref{227}, with $A>\frac1{2\rho}$. Then,  $\bar x$ solves the Euler-Lagrange equation
		$$
		2A\nabla^2F(\bar x)\nabla F(\bar x)+\nabla G(\bar x)=0.
		$$
		Multiply the previous equality by $\nabla F(\bar x)$, to conclude again, as in \eqref{218}, that $\bar x= x^*$. Again, this implies inequality~\eqref{216}. 
		
		This last proof is quite interesting since we completely avoid using the gradient flow.  Moreover, methods based on optimization problems are often robust. 
	\end{enumerate}
	%
	
	
	\subsection{Gradient flows in the space of probability measures}
	\label{sec-3.2}
	
	In this section, we reproduce the three methods of Section~\ref{sec-3.0}, this time in the space of probability measures over $M$. Before doing so, we need to introduced Otto's calculus, the main point of our method.  
	For simplicity, all computations are given on  a $d$-dimensional smooth, connected and compact Riemannian manifold $(M,\mathfrak g)$. But they can be easily generalized to the setting of weighted Riemannian manifold under the $CD(\rho,n)$ condition~\eqref{100}, as in Theorem~\ref{th:liouville} or Theorem~\ref{th:B}.

	\subsubsection{Otto's calculus}
	\label{sec-3.2.1}

	Otto's calculus, so called by C. Villani in his book~\cite{villani2009}, is a very efficient tool to compute the second derivative of a functional along its probability gradient flow.  This calculus has been developed in the seminal papers~\cite{jko1998,otto2001, otto-villani2000}. It allows to  view the space of probability measures on a manifold, at least formally, as an infinite dimensional Riemannian manifold.  Our presentation is based on~\cite{gentil-leonard2020}, to which we refer for more details (see also~\cite{gentil} for an informal presentation in French). The calculus can be viewed as a heuristic guideline but all the results can be turned into rigorous statements, see the monograph~\cite{gigli2012}.
	
	Let $\mathcal P_2(M)$ denote the space of probability measures on $M$ admitting a second moment\footnote{Since we assumed for simplicity that $M$ is compact, all probability measures admit a second moment and so  $\mathcal P_2(M)= \mathcal P(M)$ in this case.}.  
	Equip $\mathcal P_2(M)$ with the Wasserstein distance, defined as follows: for every $\mu,\nu\in\mathcal P_2(M)$, 
	$$
	W_2(\mu,\nu)=\inf \sqrt{\iint {\bf d}(x,y)^2d\pi(x,y)},
	$$
	where the infimum is taken over all transportation plans $\pi\in\mathcal P (M\times M)$ with marginals $\mu$ and $\nu$ and where $\bf d$ is the Riemannian distance of $M$. 
	
	Following the presentation of ~\cite[Chap.~1]{ambrosio-gigli2008},  a path $[0,1]\ni t\mapsto \nu_t\in\mathcal P_2(M)$ is absolutely continuous with respect to the Wasserstein distance if  
	$$
	|\dot\nu_t|:=\underset{s\rightarrow t}{\limsup}\frac{W_2(\nu_t,\nu_s)}{|t-s|}\in L^1([0,1]).
	$$ 
	It turns out that given any absolutely continuous path $(\nu_t)_{t\in[0,1]}$, there exists a unique vector field $(t,x) \mapsto V_t(x)$ in $M$, such that $\int|V_t|^2d\nu_t<\infty$ 
	and $ |\dot\nu_t|^2=\int|V_t|^2d\nu_t$ a.e. in $[0,1]$, see \cite{ambrosio-gigli2008}. In addition, the vector field $V_t$ is the limit in $L^2(\nu_t)$ of the gradient of functions  
	$\varphi_n\in C^\infty(M)$ and the continuity equation holds in the sense of distributions:
	\begin{equation}
	\label{301}
	\partial_t\nu_t+\nabla\cdot(\nu_t V_t)=0\quad\text{in $\mathcal D'(M\times(0,1))$.}
	\end{equation}
	Conversely, given any such vector field $V_t$, there exists an absolutely continuous path $(\nu_t)_{t\in[0,1]}$ such that the continuity equation \eqref{301} holds.
	In other words, for almost every $t \in [0,1]$, we may see $V_t$ as a tangent vector along the path $(\nu_t)_{t \in [0,1]}$. So, we denote 
	\begin{equation}
	\label{300}
	\dot \nu_t:=V_t
	\end{equation}
	and call $\dot \nu_t$ the velocity of the path $(\nu_t)_{t \in [0,1]}$ at time $t$. The tangent space at a point $\mu\in\mathcal P_2(M)$ can then be defined by
	\[T_{\mu} \mathcal P_2(M)=\overline{\{\nabla \varphi,\,\,\phi:M \mapsto\R,  \varphi \in C^{\infty}(M)\}}^{L^2(\mu)}\]
	and a natural Riemannian metric can be defined via the scalar product in $L^2(\mu)$ by
	$$
	\langle\nabla \varphi,\nabla\psi\rangle_\mu=\int \nabla\varphi\cdot\nabla\psi\; d\mu= \int \Gamma(\varphi,\psi)d\mu,\,\,\,\quad\text{for}\,\,\,\nabla \varphi,\nabla\psi\in T_{\mu} \mathcal P_2(M).
	$$ 
	We shall write $|\nabla \varphi|_\mu^2=\int \Gamma(\varphi)d\mu$ the corresponding Riemannian length. 
	Such a metric is often referred to as {\it the Otto metric}. In addition, thanks to the Benamou-Brenier formulation, the Wasserstein distance is the Riemannian distance associated to the Otto metric. 
	
	\subsubsection{Differentiating twice Renyi's entropy using Otto's calculus}
	
	To lighten notations and formulas,  we identify henceforth measures and densities. All the measures considered in this section are supposed to be smooth and absolutely continuous with respect to the Riemannian measure on $M$. Unless specified, all integrals are viewed with respect to the normalized Riemannian measure. 
	
	Now, we consider our main flow $(\mu_t)_{t\geq0}$, started from a probability measure $\mu_0=\mu$ and solving the following nonlinear diffusion equation
	\begin{equation}
	\label{303}
	\partial_t \mu_t = \frac{1}{\alpha}\Delta\mu_t^\alpha=\nabla\cdot \left(\mu_t \frac{1}{\alpha-1}\nabla \mu_t^{\alpha-1}\right),
	\end{equation}
	where $\alpha>0$, $\alpha\neq 1$. If the initial datum $\mu_0$ is chosen smooth, bounded and bounded away from zero, then $\mu_t$ is smooth and globally defined.\footnote{for precise statements, see Section 11.5.1 in \cite{vazquez} for the existence of a unique weak solution and the proof of Proposition 7.21 in the same book for its regularity. For a precise proof assuming only standard nonlinear parabolic regularity theory (as developped in \cite{ladizenskaia}), the interested reader can easily adapt the proof presented in Section 4.3 of \cite{zugmeyer}.}
	Then, according to the continuity equation \eqref{301}, the velocity of this flow is given by
	\begin{equation}
	\label{304}
	\dot\mu_t=-\frac{1}{\alpha-1}\nabla \mu_t^{\alpha-1}\in T_{\mu_t} \mathcal P_2(M) 
	\end{equation}
	Consider now the {\it R\'enyi entropy} (of order $\alpha>0$ with $\alpha\neq 1$),
	\begin{equation}
	\label{305}
	\mathcal R_\alpha(\mu)=\frac{1}{\alpha(\alpha-1)}\int \mu^{\alpha},\,\,\,\mu\in\mathcal P _2(M), 
	\end{equation} 
	which is the main functional used in this article. Then  the gradient of $\mathcal R_\alpha$ is given by 
	\begin{equation}
	\label{306}
	\grad_\mu \mathcal{R}_\alpha:= \frac{1}{\alpha-1}\nabla\mu ^{ \alpha-1}\in T_{\mu} \mathcal P_2(M),
	\end{equation}
	see for instance~\cite[Sec.~3.2]{gentil-leonard2020}. So, if  $(\mu_t)_{t\geq0}$ is a solution of ~\eqref{303}, then 
	$$
	\dot\mu_t=-\grad_{\mu_t} \mathcal{R}_\alpha.
	$$
	In other words,~\eqref{303} is the gradient flow of the R\'enyi entropy with respect to the Otto metric.  This was proved rigorously  in~\cite{otto2001}. 
	Furthermore, the Riemannian structure given to $ \mathcal P_2(M)$ allows us to define the covariant derivatives and the Hessian of a functional. A remarkable fact is that the Hessian of R\'enyi's entropy in the sense of Otto's calculus has an explicit formulation: for any $\mu \in \mathcal P_2(M)$ and $\nabla \phi \in T_{\mu} \mathcal{P}_2(M)$,
	\begin{equation}
	\label{307}
	\mathrm{Hess}_{\mu} \mathcal{R}_\alpha(\nabla\phi , \nabla\phi)
	=\frac{1}{\alpha}\int \SBRA{(\alpha-1)(\Delta\phi)^2+\Gamma_2(\phi)}\mu^\alpha,
	\end{equation}
	where the operator $\Gamma_2$ has been defined in~\eqref{eq:gamma22} (see~\cite{otto2001} or \cite[Sec.~3.3]{gentil-leonard2020}).
	
	\medskip

	Let us now turn to our three methods to prove inequality~\eqref{sobolev:M}, under a lower bound of the Ricci curvature.

	\subsubsection{Method based on a convex inequality for the R\'enyi entropy}
	\label{sec-3.2.2}

	We mimic the first proof proposed in Section~\ref{sec-3.0}  by using the R\'enyi entropy and the fast diffusion flow. Replace the entropy $F$ of Section \ref{sec-3.0} by $R_\alpha$, with $\alpha=1-\frac1d$ and $G$ by $-R_\beta$, with $\beta=1-\frac2d$. Then, letting $\Phi=\frac{1}{\alpha-1}\mu ^{ \alpha-1}$, it follows from \eqref{306}, \eqref{307} and the $CD(\rho,d)$ condition\footnote{Recall that $d$ is the dimension of $M$.} that
	$$
	\mathrm{Hess}_{\mu} \mathcal{R}_\alpha(\grad_{\mu} \mathcal{R}_\alpha, \grad_{\mu} \mathcal{R}_\alpha)
	=\frac{1}{\alpha}\int \SBRA{(\alpha-1)(\Delta\Phi)^2+\Gamma_2(\Phi)}\mu^\alpha\ge \frac\rho\alpha\int\Gamma(\Phi)\mu^\alpha
	$$
	while, since $\beta-3=2\alpha-4$,
	$$
	-\langle\grad_{\mu} \mathcal{R}_\alpha, \grad_{\mu} (-\mathcal{R}_\beta) \rangle_\mu= \frac1{(\alpha-1)(\beta-1)} \int \nabla\mu^{\alpha-1}\nabla\mu^{\beta-1}d\mu = \int \mu^{\alpha+\beta-3}\vert\nabla\mu\vert^2=\int \Gamma(\Phi)\mu^{\alpha}
	$$
	and so we have the exact analogue of \eqref{215}, that is. 
	\begin{equation}\label{rr}
	\mathrm{Hess}_{\mu} \mathcal{R}_\alpha(\grad_{\mu} \mathcal{R}_\alpha, \grad_{\mu} \mathcal{R}_\alpha) \ge -\frac\rho\alpha \langle\grad_{\mu} \mathcal{R}_\alpha, \grad_{\mu} (-\mathcal{R}_\beta) \rangle_\mu.
	\end{equation}
	Since $\mu^*=1$ is the unique critical point of $\mathcal{R}_\alpha$, repeating the elementary analysis\footnote{In so doing, one should restrict to a smooth, bounded, bounded away from zero initial datum $\mu_0$, so that standard nonlinear parabolic regularity theory and the maximum principle apply. In particular, the family $(\mu_t)_{t\ge0}$ is uniformly bounded and compact in the $C^k$ topologies.} leading to \eqref{214}, one has the following limits  
	\begin{equation}
	\label{309}
	\left\{
	\begin{array}{l}
	\displaystyle\lim_{t\rightarrow\infty}\mu_t=1, \\
	\displaystyle\lim_{t\rightarrow\infty}\frac{d}{dt}\mathcal R_\alpha(\mu_t)=0.
	\end{array} 
	\right.
	\end{equation}
	Hence, by the very same proof of Section \ref{sec-3.0}, we arrive at the exact analogue of \eqref{216}, that is:
	$$
	-\mathcal R_\beta(\mu^*) \le \frac\alpha{2\rho}\vert  \grad_{\mu} \mathcal{R}_\alpha\vert_\mu^2 -\mathcal R_\beta(\mu).
	$$
	By using the very definitions of $\mathcal R_\alpha$, $\mathcal R_\beta$, $\alpha$, $\beta$ and $\Phi$ we obtain 
	$$
	1\leq \frac{4(d-1)}{\rho d(d-2)}\int \Gamma(\mu^{\frac{d-2}{2d}})+\int\mu^{\frac{d-2}{d}},
	$$
	for any probability measure $\mu$. 
	Letting $\vert f\vert=\mu^{\frac{d-2}{2d}}$   in the previous inequality, we obtain 
	$$
	1\leq \frac{4(d-1)}{\rho d(d-2)}\int \Gamma(f)+\int f^2,
	$$
	under the normalization $\Vert f\Vert_{2^*}=1$ (so that $\mu$ is a probability measure). This is precisely Sobolev's inequality~\eqref{sobolev:M}.
	This proof was first proposed by J. Demange in~\cite{demange2008}. This method is important since it shows that Sobolev's inequality under a lower bound on the Ricci tensor is just a convex inequality applied to a functional (the R\'enyi entropy) along its gradient flow (the fast diffusion equation). The drawback of this method  is that it is not so easy to prove the existence of a smooth global solution of the nonlinear diffusion equation~\eqref{303} and the two limits~\eqref{309}.

	\subsubsection{Method based a minimization problem associated with the fast diffusion equation}
	\label{sec-3.3}
	
	Now, let us mimic the second proof of Section~\ref{sec-3.0}. Given $A>0$, we consider the minimization problem
	\begin{equation}\label{eqI}
	I(A):=\inf_{\mu\in\mathcal P_2(M)}\SBRA{A |\nabla \mathcal R_\alpha(\mu)|^2-\mathcal R_\beta(\mu)}
	\end{equation}
	And we prove that  for any $A>\frac\alpha{2\rho}$, $-\mathcal R_\beta(\mu^*)\le I(A)$, where $\mu^*=1$. Then, Sobolev's inequality follows as discussed in the previous section. Since the problem is critical, the first delicate point consists in proving that the infimum $I(A)$ is attained by some measure $\overline\mu$, which we admit here.\footnote{In our proof in Section 2.1, we bypassed this issue by approximating the inequality with a subcritical inequality}.This being said, once we have a well-defined global smooth solution of the gradient flow \eqref{303}, and once we've observed the strict convexity of $\mathcal R_\alpha$, which follows from \eqref{307} and the $CD(\rho,d)$ condition, then all computations done in Section \ref{sec-3.0} remain unchanged, leading to $\overline\mu=\mu^*=1$ and the desired inequality is proved.  The main advantage of this method, compared to the previous one, is that it is no longer necessary to prove the two delicate limits of the fast diffusion equation~\eqref{309}.  However, one needs to prove the existence of the minimizer $\overline\mu$ as well as the existence of a smooth solution of the fast diffusion equation~\eqref{303}.
	The method proposed in the proof of Theorem~\ref{th:sobolev} avoids both problems by working in a subcritical setting and by using the limit case, that is, the elliptic equation.
	\subsubsection{Method based only on the minimization problem}
	\label{sec-3.4}
	Indeed, mimic the third proof of Section~\ref{sec-3.0}. We consider again the minimization problem \eqref{eqI}. Assume that there exists a probability measure $\overline\mu$ minimizing $I(A)$. Then, $\overline\mu$ satisfies the corresponding Euler-Lagrange equation, given by 
	\begin{equation}
	\label{400} 
	T_{\overline\mu}\mathcal P_2(M)\ni 2A\Hess_{\overline\mu}\mathcal R_\alpha(\grad_{\overline\mu}\mathcal R_\alpha,\cdot)-\grad_{\overline\mu}\mathcal R_\beta=0,
	\end{equation}
	thanks to Otto's calculus.
	Apply the equality to the test function $\grad_{\overline\mu}\mathcal R_\alpha$, to get 
	$$
	2A\Hess_{\overline\mu}\mathcal R_\alpha(\grad_{\overline\mu}\mathcal R_\alpha,\grad_{\overline\mu}\mathcal R_\alpha)-\langle\grad_{\overline\mu}\mathcal R_\beta,\grad_{\overline\mu}\mathcal R_\alpha\rangle_{\overline\mu}=0,
	$$
	Using again the strict convexity of $\mathcal R_\alpha$ and \eqref{rr}, we conclude that $\overline\mu=1$.
	
	\medskip
	
	The proof proposed in Section~\ref{sec-2.1} is inspired from this one. The only difference is that we work here on the space of probability measures, whereas in Section~\ref{sec-2.1}, to prove the existence of a minimizer, we work on the space of functions $v$ such that $||v||_{q}=1$, where $q\in[1,2^*)$ is subcritical. The elliptic equation~\eqref{10} is, up to a change of functions, the equation~\eqref{400} whereas when we multiply  by $\Delta \Phi^{1-d'}$ and integrate in the proof of Section~\ref{sec-2.1} is exactly applying~\eqref{400} to $\grad_{\overline\mu}\mathcal R_\alpha$.

	\small
	\newcommand{\etalchar}[1]{$^{#1}$}
	
	\small
	This work was supported by the French ANR-17-CE40-0030 EFI project. 
	
	\medskip
	
	L. D., I. G. {Institut Camille Jordan, Umr Cnrs 5208, Universit\'e Claude Bernard Lyon 1, 43 boulevard du 11 novembre 1918, F-69622 Villeurbanne cedex.
		\texttt{dupaigne, gentil@math.univ-lyon1.fr}
		
		\medskip
		
		S. Z.  MAP5, UMR CNRS 8154, Universit\'e de Paris, 
		45 rue des Saints-Pères, 75270 Paris cedex 06.
		\texttt{simon.zugmeyer@u-paris.fr}

	\end{document}